\documentclass{amsart}
\usepackage{CJKutf8}
\usepackage{amsmath,amsfonts,amssymb,amscd,bm,latexsym}
\usepackage[colorlinks=true]{hyperref}
\usepackage{tikz}

\usepackage{pgflibraryarrows}
\usepackage{pgflibrarysnakes}
\usepackage{pgflibraryarrows}
\usepackage{pgflibrarysnakes}
\hypersetup{urlcolor=blue, citecolor=blue}
\numberwithin{equation}{section} 

\def\<{\langle}             \def\>{\rangle}

\newtheorem{thm}{Theorem}[section]

\newtheorem{lem}[thm]{Lemma}

\theoremstyle{definition}

\newcommand{\beeq}{\begin{equation}}\newcommand{\eneq}{\end{equation}}
\newcommand{\al}{\alpha}    \newcommand{\be}{\beta}
\newcommand{\de}{\delta}    
\newcommand{\ep}{\varepsilon}
  
    \newcommand{\la}{\lambda}
    
\newcommand{\om}{\omega}    
\newcommand{\ga}{\gamma}    
\newcommand{\R}{\mathbb{R}}

\newcommand{\Sp}{\mathbb{S}}
\newenvironment{prf}{\noindent {\bf Proof.} }{\endprf\par}
\def \endprf{\hfill  {\vrule height6pt width6pt depth0pt}\medskip}
\newcommand{\pa}{\partial}
\newcommand{\les}{{\lesssim}}
\newcommand{\supp}{\,\mathop{\!\mathrm{supp}}}

\newcommand{\gm}{\mathfrak{g}}

\numberwithin{equation}{section}
\title[Wave equation with damping and potential]
      {Lifespan estimates for wave equations with damping and potential posed on asymptotically Euclidean manifolds}
\author{Mengyun Liu}
\address{Department of Mathematics\\                	
Zhejiang Sci-Tech University\\                Hangzhou 310018, P. R. China}
\email{mengyunliu@zstu.edu.cn}

\thanks{The author was supported by NSFC 11971428.}
\date{\today}
\dedicatory{} \commby{}
\begin{document}
\begin{abstract}
In this work, we investigate the problem of finite time blow up as well as the upper bound estimates of lifespan for solutions to small-amplitude semilinear wave equations with time dependent damping and potential, and mixed nonlinearities $c_1 |u_t|^p+c_2 |u|^q$, posed on asymptotically Euclidean manifolds,
 which 
is related to both the Strauss conjecture and the Glassey conjecture. 
 \end{abstract}



\maketitle

\section{Introduction}

Let $(\R^{n}, \gm)$ be a asymptotically Euclidean (Riemannian) manifold, with $n\ge 2$. 
By asymptotically Euclidean, we mean that
$(\R^n, \gm)$ is certain perturbation of the Euclidean space $(\R^n, \gm_0)$. More precisely, we assume
$\gm$ can be decomposed as 
\begin{align}
\label{dl1}
\gm=\gm_{1}+\gm_{2}\ ,
\end{align}
where $\gm_{1}$ is a spherical symmetric, long range perturbation of $\gm_0$, and $\gm_{2}$ is an exponential (short range) perturbation.
By definition, there exists polar coordinates $(r,\omega)$ for $(\R^n,\gm_1)$, in which we can write
\beeq
\label{dl3}
\gm_{1}=K^{2}(r)dr^{2}+r^{2}d\omega^{2}\ , 
\eneq
where $d\om^2$ is the standard metric on the unit sphere $\Sp^{n-1}$, and
\beeq
\label{dl2}
|\pa^{m}_{r}(K-1)|\les  \<r\>^{-m-\rho}
,  m=0, 1, 2.
\eneq
for some given constant $\rho >0$. 
Here and in what follows, 
$\langle x\rangle=\sqrt{1+|x|^2}$, and
we use $A\les B$ ($A\gtrsim B$) to stand for $A\leq CB$ ($A\geq CB$) where the constant $C$ may change from line to line. 
Equipped with the coordinates $x=r\omega$, we have
$$\gm=g_{jk}(x)dx^j dx^k\equiv \sum^{n}_{j,k=1}g_{jk}(x)dx^j dx^k\ ,\ \gm_2=g_{2, jk}(x)dx^j dx^k\ ,$$
where we have used the convention that Latin indices $j$, $k$ range from $1$ to $n$ and the Einstein summation convention for repeated upper and lower indices. 
Concerning $\gm_2$, we assume
it is an exponential (short range) perturbation of $\gm_1$, that is, there exists $\be>0$ so that
\beeq
\label{dlfjia}
|\nabla g_{2,jk}|+|g_{2,jk}|\les e^{-\be\int^{r}_{0}K(\tau)d\tau},\ 
|\nabla^2 g_{2,jk}|\les 1\ .
\eneq
By asymptotically Euclidean and Riemannian assumption, it is clear that 
there exists a constant $\delta_{0}\in(0, 1)$ such that 
\beeq
\label{unelp}
\delta_{0}|\xi|^2\le
g^{jk}\xi_{j}\xi_{k}\le
\delta_{0}^{-1} |\xi|^2, \forall \ \xi\in\R^{n},\ 
K\in (\delta_{0}, 1/\delta_{0})
\ .
\eneq

In this paper, we are interested in the investigation of the blow up part of the wave equations with time dependent damping and potential on asymptotically Euclidean manifolds.
More precisely, we will study the blow up of solutions for 
the following semilinear wave equations with small initial data,
  posed on asymptotically Euclidean manifolds \eqref{dl1}-\eqref{unelp}
\beeq
\label{nlw}
\begin{cases}
\pa^{2}_{t}u-\Delta_{\gm}u+\frac{\mu_1}{1+t}u_t+\frac{\mu_2}{(1+t)^2}u=c_1 |u_{t}|^{p}+ c_2 |u|^q\\
u(0,x)=\ep u_0(x), u_{t}(0,x)=\ep u_1(x)\ ,
\end{cases}
\eneq
where, 
$\Delta_{\gm} =\nabla^j\partial_{j}$ is the standard Laplace-Beltrami operator,  $p, q>1$, $\mu_1, \mu_2\in \R$, $c_1, c_2\ge 0$ and $\ep>0$ is a small parameter. Concerning the initial data, we assume 
\beeq
\label{hs2}
0< u_0, u_1 \in C^{\infty}_{0}(\R^{n}),\quad \supp(u_0, u_1)\subset\{x\in \R^{n}; |x|\leq R_{0}\}\ , 
\eneq
 for some $R_{0}>0$.

 Before stating our results, let us briefly review the history of this problem. \\
 
 {\bf (\uppercase\expandafter{\romannumeral1}) No damping and potential ($\mu_1=\mu_2=0$), Euclidean space $\gm=\gm_0$}\\
 
  When $c_1=0<c_2$, the problem is related to the Strauss conjecture, and the critical power is given by the Strauss exponent $q_S(n)$, where $q_S(n)$ is the positive root of the equation:
$$(n-1)q^{2}-(n+1)q-2=0
\ .$$
Moreover, when there are no global solutions, and initial data are nontrivial and nonnegative, it has been proved that the upper bound of the lifespan is
 \beeq\label{eq-life}
T_\ep\leq S_{\ep}(n):=\left\{
\begin{array}{ll }
  C_0 \ep^{\frac{2q(q-1)}{(n-1)q^{2} - (n+1)q - 2}}    & 1<q<q_S(n);   \\
\exp(C_0\ep^{-q(q-1)})      &   q=q_S(n).
\end{array}\right.
\eneq
We refer
\cite{LW2019} for discussion of the history and related results.
 
 When $c_2=0<c_1$, the problem of determining the sharp range of powers of $p>1$ for global existence versus blow up for arbitrary small initial data, is known as the Glassey conjecture, where the critical power $p$ for \eqref{nlw}
is conjectured to be
$$p_G(n):=1 + \frac{2}{n-1}\ .$$
See\cite{FJ1981}, \cite{JK84}, \cite{sideris}, \cite{hk}, \cite{nt}, while 
nonexistence of global solutions for $1<p\le p_G(n)$ and any $n\ge 2$ as well as the upper bound of lifespan,
\beeq
\label{main-es}
T_{\ep}\leq G_{\ep}(n):=
\begin{cases}
C_{0}\ep^{-\frac{2(p-1)}{2-(n-1)(p-1)}}, \ 1<p<p_{G}(n)\ ,\\
\exp(C_{0}\ep^{-(p-1)}), \ p=p_{G}(n)\ ,
\end{cases}
\eneq
has also been well-known (at least for the case $u_1\neq 0$), see Zhou \cite{zhou} and references therein.

 When $c_1c_2\neq 0$, the problem
 is related to both the Glassey conjecture and the Strauss conjecture. It turns out that a new critical curve occurs in the expected region of global existence $p>p_G(n)$, and $q>q_S(n)$, that is
 $$\la(p, q, n)=(q-1)((n-1)p-2)=4\ , q>q_{S}(n), p>p_G(n)\ ,$$
 where the critical and super-critical case are known to admit global existence, at least for $n=2, 3$, see Hidano-Wang-Yokoyama \cite{HWY}, while there is non-existence of global existence for the sub-critical case, as well as 
an upper bound of the lifespan 
  \beeq
 \label{mixed-lifespan1}
  T^{p,q}_{\ep}\leq Z_{\ep}(n):=C_0 \ep^{-\frac{2p(q-1)}{4-\la(p, q, n)}},
  \eneq
  in the region
\beeq
\label{mixed-region1} \la(p, q, n)<4,
   p\leq\frac{2n}{n-1},  q< \frac{2n}{n-2}\ , \eneq
   when $u_1\neq 0$,
  see Han-Zhou \cite{HZ2014}. Note that,  in the work of Lai-Takamura \cite{LTmix}, they remove the restriction for $p$ and obtain an upper bound when $p>\frac{2n}{n-1}$ and
$u_1\neq 0$. 
\\
 
 {\bf (\uppercase\expandafter{\romannumeral2}) No damping and potential ($\mu_1=\mu_2=0$), Asymptotically Euclidean space}
 \\
 
 There is not much theory for the finite time blow-up phenomenon on asymptotically Euclidean space 
as well as the 
Schwarzschild/Kerr black hole spacetimes. 

When $c_1=0<c_2$, for asymptotically Euclidean manifolds
with $\gm=\gm_0+\gm_2$, Wakasa-Yordanov
 \cite{WaYo18-1pub} proved blow up results in the critical case $q=q_S(n)$. While for more general asymptotically Euclidean manifolds $\gm=\gm_1+\gm_2$, Liu-Wang \cite{LW2019} proved the blow-up results as well as upper bound of lifespan when $1<q\leq q_{S}(n)$. For the Schwarzschild black hole spacetime, Lin-Lai-Ming \cite{LinLaiMing19} obtained blow up result for $1<q\le 2$, while Catania-Georgiev \cite{CaGe06} obtained a weaker blow up result for $1<q< q_S(3)$.  See, e.g., Zha-Zhou\cite{ZZ15} for more discussion on exterior domains.

 When $c_1c_2\neq 0$ and $\gm=\gm_1+\gm_2$, 
Liu-Wang \cite{LW2020} proved that there are no global solutions to
 \eqref{nlw},
 for arbitrary small $\ep>0$, provided that
$p\le p_G(n)$, or $q\leq q_S(n)$, or 
\beeq
\nonumber
\la(p, q, n)<4\ ,
 p>1,  q >1\ .
 \eneq

 {\bf (\uppercase\expandafter{\romannumeral3}) Time dependent damping and potential $\mu_1\mu_2\neq 0$, Euclidean space} \\
 
 When there is time dependent damping and potential $\mu_1\mu_2\neq 0$, the critical exponent for \eqref{nlw} is related with two kinds of the dimensional shift due to the critical damping and potential coefficients. An important value of this problem is
 $$\de=(\mu_1-1)^2-4\mu_2\ .$$
For $c_1=0<c_2$, $\mu_1,\mu_2>0$ and $\de\in (0, 1]$, Palmieri-Reissig \cite{PR} obtained the blow-up result for 
 $$1<q\leq q_{cri}=\max\{q_F(n+\al)\ , q_{S}(n+\mu_1)\}\ , \al=\frac{\mu_1-1-\sqrt{\de}}{2}\ ,$$
except for $n=1$ and $q=q_S(n+\mu_1)$, where $q_F(n)$ is Fujita exponent (see, e.g., \cite{LZ, Fujita66}) satisfies
$$\gamma_F(q, n): =2-n(q-1)=0\ ,$$
for heat or damped wave equations
\[
u_{tt}+u_t-\Delta u=|u|^p, u_t-\Delta u=|u|^p\ .
\] 
Later, Palmieri-Tu \cite{PT} relax the restriction of $\de$ and showed the blow-up results for $1<q\leq q_S(n+\mu_1)$ with $\mu_1, \mu_2, \de \geq 0$. Recently, Lai-Schiavone-Takamura \cite{LST} remove the nonnegative restriction on $\mu_1, \mu_2$ and obtained the blow-up results for $1<q<q_{cri}$ with $\de\geq 0$ and $\mu_1, \mu_2\in \R$. For $c_2=0<c_1$ and $\mu_1, \mu_2\geq 0$, Palmieri-Tu \cite{PT2} obtained the blow-up results for $1<p\leq p_G(n+\sigma)$, where 
\begin{equation*}
\sigma=
\left\{
\begin{array}{ll}
\mu_1+1-\sqrt{\de}, & \mathrm{if}~\de\in[0, 1)\ ,\\
\mu_1,& \mathrm{if}~\de\geq 1\ .
\end{array}
\right.
\end{equation*}
Recently, Hamouda-Hamza \cite{HH} made the improvement by extending the blow-up region to $1<p\leq p_G(n+\mu_1)$ for any $\de\geq 0$.  For the case $\mu_1\mu_2=0$, see, e.g., \cite{LST}, \cite{PT2}, \cite{LTW} and references therein for the more discussion on the history. 

While for $c_1c_2\neq 0$, Hamouda-Hamza \cite{HH} proved the blow-up results for 
$$\la(p, q, n+\mu_1)<4,\  p>p_G(n+\mu_1), \ q>q_S(n+\mu_1)\ ,$$
with $\mu_1, \mu_2, \de\geq 0$ and the upper bound lifespan estimate $$T_{\ep}\  \les \ Z_{\ep}(n+\mu_1)\ .$$ 
 
The main result of this paper then states that there are no global solutions to
 \eqref{nlw}
 for arbitrary small $\ep>0$, provided that
$p\le p_G(n+\mu_1)$, or $q<\max\{q_F(n+(\mu_1-1-\sqrt{\de})/2)\ , q_{S}(n+\mu_1)\} $, or 
\beeq\label{mixed}
\la(p, q, n+\mu_1)<4,\ 
 p>1, \  q >1\ .\eneq
More precisely, we have
\begin{thm}
\label{main}
Let $n\geq 2$, $p, q>1$, $\mu_1, \mu_2\in \R$  and $\de\geq 0$. Consider \eqref{nlw} with $c_1, c_{2}\ge 0$ and nontrivial initial data \eqref{hs2} and 
$$\Big(\al\int u_0 dv_{\gm}+\ep\int u_1 dv_{\gm}\Big)\geq 0\ ,$$
posed on asymptotically Euclidean manifolds \eqref{dl1}-\eqref{unelp}. Suppose it has a weak solution $u\in C^{2}([0, T_\ep);\mathcal{D}'(\R^{n}))$ with $u_{t}\in L_{loc}^{p}([0, T_{\ep})\times \R^{n})$, $|u|^{q}\in C([0, T_\ep);\mathcal{D}'(\R^{n}))$ and
\beeq
\label{fsp}
\supp u\subset \{(t, x); \int^{r}_{0}K(\tau)d\tau\leq t+R_{1}\},
\eneq
for some $R_{1}\ge \int_{0}^{R_{0}}K(\tau)d\tau$. Then 
there exist constants $\ep_0>0$ and $C_{0}>0$, such that for any $\ep\in (0, \ep_{0})$, we have the following results on the upper bound:
\begin{equation*}
T_{\ep}\leq 
\left\{
\begin{array}{ll}
G_{\ep}(n+\mu_1), & \mathrm{if}~c_{1}>0,  1<p\leq p_{G}(n+\mu_1)\ ,\\
S_{\ep}(n+\mu_1),& \mathrm{if}~ c_{2}>0,  1<q< q_{S}(n+\mu_1)\ ,\\
 C_0\ep^{-\frac{q-1}{\ga_F(q, n+\al)}}, &\mathrm{if}~ c_{2}>0,  1<q< q_F(n+\al)\ ,\\
Z_{\ep}(n+\mu_1),& \mathrm{if}~c_1c_2>0\ , \la(p, q, n+\mu_1)<4\ , p, q>1\ .\\
\end{array}
\right.
\end{equation*}

\end{thm}

 \subsubsection*{Outline} Our paper is organized as follows. In Section \ref{sec:1stTest}, we collect the Kato type lemma and the special solutions for the elliptic ``eigenvalue" problems \eqref{1.1}, with certain asymptotic behavior. These solutions play a key role in constructing the test functions and the proof of blow up results. In Section \ref{sec:AE3}, we construct the test function and collect their properties.  Then, in Section \ref{proof}, we give the proof of Theorem \ref{main}, by applying a relatively routine argument (see, e.g.,
 \cite{LW2020} or \cite{HH}).

  \section{Preliminary}\label{sec:1stTest} 
  
As is standard, when we employ the test function method, we typically need  the Kato type lemma to conclude nonexistence of global solutions as well as the upper bound of the lifespan, for which proof, we refer Sideris \cite{Sideris84} for the blow up result and Zhou-Han \cite{ZhouHan11} for the upper bound. 
\begin{lem}[Kato type lemma]
\label{lem4.1}
Let $\be>1$, $a\geq 1$ and $(\be-1)a>\al -2$. If $F\in C^{2}([0, T))$ satisfies
$$F(t)\geq \delta (t+1)^{a},\ F''(t)\geq k(t+1)^{-\al}F^{\be}\ ,$$
for some positive constants $\delta, k$. Then $F(t)$ will blow up at finite time, that is, $T<\infty$. Moreover, we have the following upper bound of $T$,
$$T\leq c\delta^{-\frac{\be-1}{(\be-1)a-\al+2}},$$
for some constant $c$ which is independent of $\delta$. 
\end{lem}

 Next, we collect the asymptotic behavior of solutions for elliptic equation
 \beeq
\label{1.1}
\Delta_{\gm}\phi_\la=\lambda^{2}\phi_\la\ .
\eneq
These solutions  will play a key role in the construction of the test functions. 
  \begin{lem}[Lemma 3.1 in \cite{LW2019}]\label{elp}
Let $n\geq 2$ and 
$(\R^n, \gm)$ be asymptotically Euclidean manifold with \eqref{dl1}-\eqref{unelp}.
Then there exist $\la_0, c_0>0$ such that for any
$0<\la\le \la_0$, there is a solution of \eqref{1.1} satisfying
\beeq
\label{1.50}
c_{0} 
<\phi_\la(x) < c_0^{-1}\< \la r\>^{-\frac{n-1}{2}}e^{\la \int^{r}_{0}K(\tau)d\tau} \ .
\eneq
\end{lem}

\subsection{Finite speed of propagation}   
Since we suppose the support of  solution $u$ of \eqref{nlw} satisfies 
\beeq
\nonumber
\supp u\subset \{(t, x); \int^{|x|}_{0}K(\tau)d\tau\leq t+R_{1}\}=D\ .
\eneq
Recall that $K(r)\in [\delta_{0}, 1/\delta_{0}]$, then we have 
\beeq
\nonumber
D\subset\{(t,x); |x|\leq \frac{t}{\delta_{0}}+\frac{R_{1}}{\delta_{0}}\}=\tilde{D}\ .
\eneq

\section{Test functions}\label{sec:AE3}
In this section, we construct the test function and collect the properties of two auxiliary functions, which we shall use later.
\subsection{Solution of linear dual equation}
Let $\psi(t, x)=\rho(t)\phi_{\la_1}(x)$, where $\phi_{\la_1}(x)$ is the solution of elliptic solution of \eqref{1.1} for some constant $\la_1\in (0, \la_0)$. And $\rho(t)$ is the solution of the following ordinary differential equation
$$\rho''-(\frac{\mu_1}{1+t}\rho)'+\frac{\mu_2}{(1+t)^2}\rho=\la_1\rho\ ,$$
which can be read as 
$$\rho(t)=(1+t)^{\frac{\mu_1+1}{2}}K_{\frac{\sqrt{\de}}{2}}(\la_1(t+1))\ ,$$
here, $K_{\nu}(t)$ is the modified Bessel function of the second kind. Then it is easy to check that $\psi(t, x)$ is the solution to 
$$\pa^{2}_t\psi-\Delta_{\gm}\psi-\frac{\pa}{\pa_t}\big(\frac{\mu_1}{1+t}\psi\big)+\frac{\mu_2}{(1+t)^{2}}\psi=0\ .$$
Moreover, we also need the $L^{m}$ estimate of $\psi(t, x)$ in the domain $D$.
\begin{lem}
Let $m>1$, then we have 
\beeq
\int_{D}\psi^{m} \ dv_{\gm} \ \les \ (t+1)^{n-1-\frac{n-1-\mu_1}{2}m}\ , t\geq T_{*}\ ,
\eneq
for some $T_{*}=T_{*}(\mu_1, \mu_2)>1 $. \end{lem}
\begin{prf}
By \cite{GRE}, we have
 $$K_{\frac{\sqrt{\de}}{2}}(\la_1(t+1))\les (t+1)^{-1/2}e^{-\la_1 t}\ , t\geq T_{*}\ ,$$
 for some $T_{*}$ depends on $\de$ thus depends on $\mu_1, \mu_2$. Then by above and  \eqref{1.50} we obtain that for $t\geq T_*$
\begin{align*}
\int_{D}\psi^{m} dv_{\gm} \les (1+t)^{\frac{m\mu_1}{2}}\int_{\int^{|x|}_{0}K(\tau)d\tau\leq t+R_{1}}e^{m\la_1(\int^{|x|}_{0}K(\tau)d\tau-t)} (1+|x|)^{-\frac{n-1}{2}m} dv_{\gm}\ .
\end{align*}
We divide the region $D$ into two disjoint parts: $D=D_{1}\cup D_{2}$ where
$$D_{1}=\{(t, x); \int^{|x|}_{0}K(\tau)d\tau\leq \frac{t+R_{1}}{2}\}\ .$$
For the region $D_{1}$,  we have 
$$\int_{D_{1}}\psi^{m} d v_{\gm}\les (1+t)^{\frac{m\mu_1}{2}}e^{-\la_1 mt}\int_{D_{1}}(1+|x|)^{-\frac{n-1}{2}m}e^{m\la_1\int^{|x|}_{0}K(\tau)d\tau}d v_{\gm}\ .$$
Let $\tilde{r}=\int^{|x|}_{0}K(\tau)d\tau$, then $d\tilde{r}=K(r)dr$ and $\delta_{0}r\leq \tilde{r}\leq r/\delta_{0}$ since $K\in[\delta_{0}, 1/\delta_{0}]$. Then we get 
\begin{align*}
\int_{D_{1}}\psi^{m} d v_{\gm}&\les (1+t)^{\frac{m\mu_1}{2}}e^{-\la_1 mt}\int^{\frac{t+R_{1}}{2}}_{0}(1+\tilde{r})^{n-1-\frac{n-1}{2}m}e^{m\la_1\tilde{r}}d\tilde{r}\\
&\les (1+t)^{\frac{m\mu_1}{2}}e^{-\la_1 mt}\int^{\frac{t+R_{1}}{2}}_{0} e^{\frac{3}{2}m\la_1\tilde{r}}d\tilde{r}\\
&\les (1+t)^{\frac{m\mu_1}{2}}e^{-\frac{m\la_1}{4}t} \les (1+t)^{n-1-\frac{n-1}{2}m+\frac{m\mu_1}{2}},
\end{align*}
where we have used the fact that $e^{-t}$ decays faster than any polynomial. 
For the region $D_{2}$, it is easy to see
\begin{align*}
\int_{D_{2}}\psi^{m} d v_{\gm}&\les (1+t)^{\frac{m\mu_1}{2}}e^{-\la_1 mt}\int^{t+R_{1}}_{\frac{t+R_{1}}{2}}e^{m\la_1\tilde{r}}(1+\tilde{r})^{n-1-\frac{n-1}{2}m}d\tilde{r}\\
&\les (1+t)^{n-1-\frac{n-1}{2}m+\frac{m\mu_1}{2}}\int^{t+R_{1}}_{\frac{t+R_{1}}{2}}e^{m\la_1\tilde{r}-\la_1 mt}d\tilde{r}\\
&\les (1+t)^{n-1-\frac{n-1}{2}m+\frac{m\mu_1}{2}}\ ,
\end{align*}
which completes the proof.

\end{prf}
\subsection{Two auxiliary functions}  
Let $$G_1(t)=\int u(t, x) dv_{\gm}\ , G_2(t)=\int u_t(t, x) d v_{\gm}\ .$$
Then by \cite[Lemma 3.2, Lemma 3.4]{HH}, we have 
\beeq
\label{yong1}
G_1(t) \geq C_1\ep, \ \forall \ t\geq T_0\ ,
\eneq
for some $T_0=T_0(\mu_1, \mu_2)>0$ and there exists a small $\ep_0$, such that for any $\ep\in (0, \ep_0)$, we have 
 \beeq
 \label{yong2}
 G_2(t)\geq C_2 \ep, \ \forall\ t\geq T_1=C_{3}\ln \ep^{-1}\ .
 \eneq
Here, the constant $C_1, C_2, C_3$ depends on $u_0, u_1, n, \mu_1, \mu_2$.

 \section{Proof of Theorem \ref{main}}\label{proof}
 In this section, we give the proof of Theorem \ref{main} for the cases $c_2>0$ and $c_1c_2\neq 0$, by
 applying Kato type lemma. While for the case $c_2=0<c_1$, the proof follows \cite[section 5]{HH}.

\subsection{Finite time blow up of \eqref{nlw} with $c_2>0$}

In this section, we present a proof of \eqref{mixed-lifespan1},
for \eqref{nlw} with $c_2>0$.
For that purpose, we introduce the function
$$F(t)=\int u(t, x) d v_{\gm}\ ,$$
which, as we shall see, will blow up in finite time. 


In view of \eqref{nlw}, we have 
\beeq
\label{P-11}
F''+\frac{\mu_1}{1+t}F'+\frac{\mu_2}{(1+t)^2}F=\int \big(c_1|u_{t}|^{p}+c_2 |u|^{q}\big) d v_{\gm}=N(u)\ .
\eneq
Concerning the relation between $F(t)$ and $N(u)$, by H\"older's inequality, we have 
$$|F|\leq \int_{\tilde{D}}|u| dv_{\gm}\les \Big(\int_{\tilde{D}}|u|^{q} dv_{\gm}\Big)^{1/q}(1+t)^{n/q'},$$
thus we get that 
\beeq
\label{mixn}
N(u)\geq \int_{\tilde{D}}c_2|u|^{q} dv_{\gm}\gtrsim \frac{|F|^{q}}{(t+R)^{n(q-1)}}\ .
\eneq
Let $H(t)=F(t)(1+t)^{\al}$, then we have $H$ satisfies 
\beeq 
\label{bian1}
H''+\frac{1+\sqrt{\de}}{1+t}H'=(1+t)^{\al}N\ ,
\eneq
that is 
$$((1+t)^{1+\sqrt{\de}}H')'=(1+t)^{1+\sqrt{\de}+\al}N\geq 0\ ,$$
thus we get that 
\begin{align*}
H'(t)\geq (1+t)^{-1-\sqrt{\de}}H'(0)=&(1+t)^{-1-\sqrt{\de}}(\al F(0)+F'(0))\\
=&(1+t)^{-1-\sqrt{\de}}\Big(\ep\al\int u_0 dv_{\gm}+\ep\int u_1 dv_{\gm}\Big)\geq 0\ ,
\end{align*}
which yields that 
\beeq
\label{bian1chu}
H(t)\geq H(0)=\ep\int u_0 dv_{\gm}> 0\ .
\eneq

(\uppercase\expandafter{\romannumeral1}) $\de\geq 1$: let $L(t)=(1+t)^{\frac{1+\sqrt{\de}}{2}}H(t)$, by \eqref{bian1} we have 
 $$L'(t)=\frac{1+\sqrt{\de}}{2}(1+t)^{\frac{1+\sqrt{\de}}{2}-1}H(t)+(1+t)^{\frac{1+\sqrt{\de}}{2}}H'(t)\geq 0\ ,$$
and 
$$
L''\geq L''+\frac{1-\de}{4}\frac{1}{(1+t)^2}L=(1+t)^{\al+\frac{1+\sqrt{\de}}{2}}N\ .
$$
Concerning the relation between $L(t)$ and $N(u)$, by \eqref{mixn} we have
$$N\gtrsim \frac{H^{q}}{(t+R)^{n(q-1)+\al q}}\gtrsim \frac{L^q}{(t+R)^{(n+\mu_1/2)(q-1)}}\ .$$
Combine the initial data \eqref{bian1chu}, we get that 
\begin{align}
\label{gui1}
\begin{cases}
L''\gtrsim (1+t)^{\al+\frac{1+\sqrt{\de}}{2}}N\gtrsim\frac{L^q}{(t+R)^{(n+\mu_1/2)(q-1)}}\ ,\\
L(t)\geq L(0)\geq \ep C_{u_0} (t+1)^{(1+\sqrt{\de})/2} \ .
\end{cases}
\end{align}

(\uppercase\expandafter{\romannumeral2}) $\de< 1$: by \eqref{bian1}, we have that 
$$H''+\frac{2}{1+t}H'\geq H''+\frac{1+\sqrt{\de}}{1+t}H'=(1+t)^{\al}N\ .$$
Let $L(t)=(1+t)H(t)$, by above, we have that $L$ satisfies
\begin{align}
\label{gui2}
\begin{cases}
L''\geq (1+t)^{\al+1}N \gtrsim \frac{L^q}{(t+R)^{(n+\al+1)(q-1)}}\ ,\\
 L(t)\geq L(0)\geq \ep C_{u_0} (1+t)\ .
 \end{cases}
\end{align}
By applying Kato type Lemma \ref{lem4.1} to \eqref{gui1} or \eqref{gui2}, we obtain that the lifespan estimates of $L(t)$ 
$$T_\ep \ \les\  \ep^{-\frac{q-1}{2-(q-1)(n+\al)}} , 1<q<1+\frac{2}{n+\al}\ ,$$

To show blow up of $1<q<q_S(n+\mu_1)$, we 
need to
get another lower bound of $L(t)$, for which purpose, we use the auxiliary function $G_1(t)$. By \eqref{yong1} and H\"older's inequality, for $t\geq T_0+T_{*}$, we have 
$$\ep\  \les \ G_1 \les\  \left(\int_{D} |u|^{q} dv_{\gm}\right)^{\frac{1}{q}}\left(\int_{D}\psi^{q'}dv_{\gm}\right)^{\frac{1}{q'}}\ ,$$
thus 
$$N\geq c_2\int_{D} |u|^{q}dv_{\gm} \gtrsim \ep^q  \left(\int_{D}\psi^{q'}dv_{\gm}\right)^{1-q}\gtrsim \ep^q
(t+1)^{n-1-\frac{n-1+\mu_1}{2}q}\ .$$
Then by \eqref{gui2} (or \eqref{gui1}), we have 
$$L''\gtrsim \ep^q(1+t)^{n+\al-\frac{n-1+\mu_1}{2}q}\ .$$
Notice that, if $n+\al-\frac{n-1+\mu_1}{2}q\geq -1$, then heuristically, we can integrate the above expression twice and obtain
\beeq
L \gtrsim \ep^q(1+t)^{n+\al+2-\frac{n-1+\mu_1}{2}q}\ ,
\eneq
otherwise, by \eqref{gui2}, we still have 
$$L \gtrsim \ep(1+t)\gtrsim \ep^q(1+t)^{n+2-\frac{nq+\sqrt{\de}}{2}}\ ,$$
which could be improved with the help of \eqref{gui2}, if 
$$\Big(n+2-\frac{nq+\sqrt{\de}}{2}\Big)q-(n+\al+1)(q-1)+2>n+2-\frac{nq+\sqrt{\de}}{2}\ ,$$
that is $$q<q_S(n+\mu_1)\ .$$
Then the classical Kato type lemma 
could be applied to prove blow up results and 
the desired upper bound of lifespan 
$$T_0+T_{*}\ \les \ T_{\ep} \ \les \ S_{\ep}(n+\mu_1)\ .$$
 
\subsection{Finite time blow up of \eqref{nlw} with $c_1c_2\neq 0$}

To show blow up of $\la(p, q, n+\mu_1)<4$, we 
need to
get another lower bound of $L(t)$, for which purpose, we use the auxiliary function $G_2(t)$. By \eqref{yong2} and H\"older's inequality, for $t\geq T_1+T_{*}$, we have 
$$\ep \ \les\  G_2(t) \ \les \left(\int_{D} |u_t|^{p} dv_{\gm}\right)^{\frac{1}{p}}\left(\int_{D}\psi^{p'}dv_{\gm}\right)^{\frac{1}{p'}}\ ,$$
thus 
$$N\geq c_1\int_{D} |u_t|^{p}dv_{\gm} \gtrsim \ep^p  \left(\int_{D}\psi^{p'}dv_{\gm}\right)^{1-p}\gtrsim \ep^p
(t+1)^{n-1-\frac{n-1+\mu_1}{2}p}\ .$$
Then by \eqref{gui2} (or \eqref{gui1}) and the same argument above, we can integrate twice and obtain
\beeq
L \gtrsim \ep^p(1+t)^{n+\al+2-\frac{n-1+\mu_1}{2}p}\ ,
\eneq
which could be improved with the help of \eqref{gui2}, if 
$$\Big(n+\al+2-\frac{n-1+\mu_1}{2}p\Big)q-(n+\al+1)(q-1)+2>n+\al+2-\frac{n-1+\mu_1}{2}p\ ,$$
that is $$\la(p, q, n+\mu_1)<4\ , p>1, q>1 \ .$$
Then the classical Kato type lemma 
could be applied to prove blow up results and 
the desired upper bound of lifespan 
$$T_1+T_{*}\ \les \ T_{\ep} \ \les \ Z_{\ep}(n+\mu_1)\ .$$

\bibliographystyle{plain1}

\end{document}